# A GENERAL LOWER BOUND FOR MIXING OF SINGLE-SITE DYNAMICS ON GRAPHS

By Thomas P. Hayes[1] and Alistair Sinclair[2]

*Toyota Technological Institute and University of California at Berkeley*

We prove that any Markov chain that performs local, reversible updates on randomly chosen vertices of a bounded-degree graph necessarily has mixing time at least $\Omega(n \log n)$, where $n$ is the number of vertices. Our bound applies to the so-called "Glauber dynamics" that has been used extensively in algorithms for the Ising model, independent sets, graph colorings and other structures in computer science and statistical physics, and demonstrates that many of these algorithms are optimal up to constant factors within their class. Previously, no superlinear lower bound was known for this class of algorithms. Though widely conjectured, such a bound had been proved previously only in very restricted circumstances, such as for the empty graph and the path. We also show that the assumption of bounded degree is necessary by giving a family of dynamics on graphs of unbounded degree with mixing time $O(n)$.

**1. Introduction.** A large fraction of Markov chain Monte Carlo algorithms, as studied in both theoretical and practical settings, fall into the class of "reversible single-site dynamics," or "Glauber dynamics," on bounded-degree graphs. In this class of algorithms, one is given a finite undirected graph $G = (V, E)$ of maximum degree $\Delta$ and a finite set of values $Q$, together with a Markov chain whose state space is $\Omega \subseteq Q^V$, a subset of all possible assignments of values to the vertices of $G$. (Following the terminology of spin systems in statistical physics, which is one of the major application areas, we shall refer to vertices as "sites," values as "spins" and assignments as "configurations." Note that we do not assume that all configurations are

Received October 2005; revised November 2006.
[1]Supported by an NSF Postdoctoral Fellowship and by NSF Grant CCR-0121555.
[2]Supported in part by NSF Grant CCR-0121555.
*AMS 2000 subject classifications.* Primary 60J10; secondary 60K35, 68W20, 68W25, 82C20.
*Key words and phrases.* Glauber dynamics, mixing time, spin systems, Markov random fields.

allowed, i.e., the spin system may have so-called "hard constraints.") The Markov chain is designed to sample from a given probability distribution $\pi$ on $\Omega$. At each step, it selects a site $v \in V$ uniformly at random and modifies the spin at $v$ according to a randomized update rule; this rule is required to be *local*, in the sense that it depends only on the current spins at $v$ and its neighbors in $G$, and *reversible* with respect to $\pi$. (See Example 5 in Section 2 for precise definitions.) Under mild additional assumptions, these conditions ensure that the Markov chain is ergodic and converges to the stationary distribution $\pi$. Thus, an algorithm which starts from an arbitrary initial configuration and performs sufficiently many random spin updates samples configurations (approximately) from $\pi$.

While our framework is rather more general, the most common scenario is when $\pi$ is a *Markov random field*; that is, for any subset of sites $U \subseteq V$, when the spins outside $U$ are fixed, the conditional distribution over configurations inside $U$ depends only on the spins on the boundary, $\partial U$. Markov random fields are ubiquitous in statistical physics (where they describe systems of particles with nearest-neighbor interactions) and in statistical modeling (where they capture conditional independence relations among random variables). For Markov random fields there is a particularly natural update rule called the "heat-bath" rule (or "Gibbs sampler"), in which the spin at $v$ is replaced by a random spin chosen from the correct distribution conditional on the spins of its neighbors. Many other variants (such as the "Metropolis rule"—see Section 2) are possible.

Clearly, the efficiency of the Markov chain Monte Carlo algorithm depends on the rate of convergence to equilibrium, or "mixing time" of the Glauber dynamics. This is the number of steps required until the distribution of the Markov chain is close (in total variation distance) to the stationary distribution $\pi$. In recent years, motivated by applications in statistical physics, combinatorics and artificial intelligence, much progress has been made in deriving upper bounds on the mixing time. In many cases, it has been possible to show that the mixing time is as fast as $O(n \log n)$, where $n$ is the number of sites in $G$. Celebrated examples include the two-dimensional Ising model above the critical temperature [14], graph colorings with sufficiently many colors [13] and the hard-core model (independent sets) at sufficiently low densities [17]. An upper bound of $O(n \log n)$ arises naturally from various techniques for bounding the mixing time, such as coupling [2] and the log-Sobolev constant [7] and, at least among physicists, is generally taken as the criterion for "rapid mixing" (rather than the weaker notion of being bounded by a polynomial in $n$, which is more common in computer science).

The generally accepted "folklore" about these results is that they are *optimal*, in the sense that the mixing time of Glauber dynamics can never be $o(n \log n)$. The justification for this folklore is that, by a standard coupon collecting argument, after substantially less than $n \ln n$ steps a significant number of sites have never been updated from their initial values, so mixing cannot have occurred. However, converting this intuition into a rigorous proof has turned out to be elusive except in very restricted cases, namely when $G$ is the empty graph (i.e., all spins are completely independent) [2, 6, 8], and when $G$ is a path and $\pi$ the uniform distribution over proper $q$-colorings ($q \geq 3$) [9]. This is actually a manifestation of a relative lack of tools for proving *lower bounds* on the mixing time of Markov chains, in contrast to the wide range of techniques for upper bounds that is now available. (The standard tool is to identify a "bottleneck" in the chain, and has been used to obtain exponential lower bounds on mixing times in specific cases; see, e.g., [5, 10, 12].)

To get a feel for why this issue is more subtle than it may seem at first glance, consider the very simple case in which $G$ is the empty graph on $n$ vertices, the spins are $Q = \{0, 1\}$, and the update rule replaces the spin at the chosen site by the outcome of a fair coin toss. Thus in the stationary distribution $\pi$ all spins are i.i.d. uniform. By coupon collecting, the number of steps needed to update all the sites is, with high probability, $n \ln n$ (plus lower-order terms). However, it is well known [2] that only about $\frac{1}{2} n \ln n$ steps suffice to mix in this case, so the coupon collecting analogy fails. The intuition behind this bound is that after $\frac{1}{2} n \ln n + \omega(n)$ steps only $o(n^{1/2})$ of the spins will not have been updated, which is of lower order than the standard deviation of the number of 1-spins in the stationary distribution. Hence the effect of the untouched spins is statistically undetectable in total variation. A further twist on this example yields a more surprising outcome. Suppose we modify the update rule so that the chosen spin is flipped with probability $1 - \frac{1}{n+1}$ and left unchanged with probability $\frac{1}{n+1}$. Remarkably, it turns out [2, 6] that the mixing time for this modified dynamics is only $\frac{1}{4} n \ln n + O(n)$, after which time $\Theta(n^{3/4})$ of the sites have not yet been updated! The reason for this discrepancy, of course, is that the update distribution is not uniform: by judiciously choosing the probability of flipping the spin, we have accelerated mixing. In this case the speedup is only by a constant factor. But, as this simple example indicates, it is quite plausible that in a more complex setting with nontrivial spin interactions it may be possible to tune the Glauber dynamics so as to achieve $o(n \log n)$ mixing time.

In this paper we prove, in a very general setting, that this is *not* possible. Specifically we show that *any* nontrivial Glauber dynamics on *any* graph $G$ of bounded degree must have mixing time $\Omega(n \log n)$. (The constant concealed in the $\Omega$ here depends only on the maximum degree $\Delta$ of $G$.)

THEOREM 1.1. *Let $\Delta \geq 2$ be a fixed positive integer, and let $G$ be any graph on $n$ vertices of maximum degree at most $\Delta$. Any nonredundant[4] Glauber dynamics on $G$ has mixing time $\Omega(n \log n)$, where the constant in the $\Omega(\cdot)$ depends only on $\Delta$.*

We note that Theorem 1.1 applies, in particular, to the standard heat-bath dynamics for all widely studied nearest-neighbor systems in statistical physics, including the Ising model, the hard-core model, and graph colorings, with arbitrary boundary conditions.

Our lower bound actually holds in somewhat greater generality. In particular, it applies to so-called "block dynamics," a variant of Glauber dynamics in which at each step a randomly chosen constant-size block of sites (rather than just a single site) is updated, again via a local, reversible rule; moreover, the assumption that the site to be updated is chosen uniformly can be relaxed to the requirement that no site is chosen with probability greater than $O(1/n)$ (no lower bound on the site selection probabilities is needed). On the other hand, the assumption of bounded degree *is* essential: we complement Theorem 1.1 by exhibiting a Glauber dynamics on a family of $n$-vertex graphs of unbounded degree for which the mixing time is only $O(n)$. In fact, we prove the following trade-off between maximum degree and the mixing time:

THEOREM 1.2. *For each $n$, let $\Delta(n)$ be any natural number satisfying $2 \leq \Delta(n) < n$. Then there exists a family of graphs $G = G(n)$, where $G(n)$ has $n$ vertices and maximum degree $\Delta(n)$, and an associated Glauber dynamics on $G(n)$ with mixing time $O(n \log n / \log \Delta(n))$.*

The two main ingredients in the proof of Theorem 1.1 are "disagreement percolation," which bounds the rate at which information can flow between sites in the Glauber dynamics as a function of their distance, and "complete monotonicity," which implies that, for a given site, under suitable initial conditions, the probability of the site having its initial spin decreases monotonically with time. We use these tools, together with a two-stage coupling argument, to identify initial conditions for the dynamics which have a statistically observable effect that persists for $\Omega(n \log n)$ steps. A similar overall strategy was used by Dyer, Goldberg and Jerrum in their study of the Glauber dynamics for colorings of a path [9]. However, in their case the simplicity of the graph and the restriction to colorings made it possible to do the calculations explicitly, whereas in our general setting we have to replace those calculations with general arguments based on percolation and monotonicity in arbitrary graphs with arbitrary spin systems.

---

[4]The nonredundancy condition merely requires that for every site there are at least two distinct spin values that appear in some feasible configuration; see Section 2.1. Without this condition, we could simply pad the graph with an arbitrary number of *"frozen"* sites without increasing the mixing time per site.

The remainder of the paper is organized as follows. In the next section we introduce some basic terminology and background, including a translation from discrete- to continuous-time Glauber dynamics which allows us to use the latter for convenience in the rest of the paper. In Section 3 we develop our principal tools: disagreement percolation and complete monotonicity; and in Section 4 we use them to prove our main theorem, a lower bound of $\Omega(n \log n)$ for the Glauber dynamics on any bounded-degree graph. Finally, in Section 5 we discuss various extensions of our main result, and give examples which show that the assumptions of bounded degree and of (moderately) uniform random site selection are necessary.

## 2. Preliminaries.

2.1. *Glauber dynamics.* Let $G = (V, E)$ be an $n$-vertex undirected graph of maximum degree at most $\Delta \geq 2$, and let $Q = \{1, \ldots, q\}$ be a finite set of *spins*. We refer to the vertices of $G$ as *sites*, and we identify (induced) subgraphs of $G$ with their sets of sites. A *configuration* on $G$ is an assignment $\sigma : V \to Q$ of spins to sites. We specify an arbitrary set $\Omega \subseteq Q^V$ of configurations, which we call *feasible*. If not all configurations are feasible, we say that the associated spin system has "hard constraints."

We will use the term "Glauber dynamics" to refer to any ergodic Markov chain on $\Omega$ that, at each step, picks a site $v \in V$ uniformly at random and applies a local, reversible update to the spin at $v$. More precisely, in (feasible) configuration $\sigma \in \Omega$ it does the following:

1. pick $v \in V$ uniformly at random (u.a.r.);
2. replace the spin $\sigma(v) = s$ by a new spin $s'$ chosen with probability $\kappa_{\sigma,v}(s, s')$, resulting in a new feasible configuration $\sigma'$.

The update distributions $\kappa_{\sigma,v}(s, \cdot)$ are required to be *local*, that is, to depend only on the values of $\sigma$ at $v$ and its neighbors, and *reversible* with respect to some probability distribution $\pi$ that assigns nonzero weight to every $\sigma \in \Omega$, that is, to satisfy the detailed balance conditions

$$\pi(\sigma)\kappa_{\sigma,v}(s, s') = \pi(\sigma')\kappa_{\sigma',v}(s', s).$$

The most common scenario in applications is when the distribution $\pi$ is a "Markov random field" (MRF) on $G$. This means that, for all subsets $U \subseteq V$, when the configuration outside $U$ is fixed, then the conditional distribution on the configuration inside $U$ depends only on the spins on the *boundary* of $U$, that is, on $\partial U = \{v \in V \setminus U : \exists u \in U \text{ with } \{u, v\} \in E\}$. Thus in particular the conditional distribution $\pi_{\sigma,v}(\cdot)$ of the spin at $v$ depends only on the spins $\sigma(u)$ at the neighbors $u$ of $v$. This admits the following two very natural local update rules that are reversible with respect to $\pi$:

1. *Heat-bath dynamics* (*or Gibbs sampler*). Set $\kappa_{\sigma,v}(s,s') = \pi_{\sigma,v}(s')$, that is, the new spin at $v$ is chosen from the conditional distribution given the spins at the neighbors of $v$.

2. *Metropolis dynamics.* Set $\kappa_{\sigma,v}(s,s') = \frac{1}{|Q|} \min\{\frac{\pi_{\sigma,v}(s')}{\pi_{\sigma,v}(s)}, 1\}$ for $s \neq s'$, and $\kappa_{\sigma,v}(s,s) = 1 - \sum_{s' \neq s} \kappa_{\sigma,v}(s,s')$. That is, choose $s' \in Q$ u.a.r., and "accept" this choice with probability that depends on the ratio $\pi_{\sigma,v}(s')/\pi_{\sigma,v}(s)$, otherwise leave $\sigma(v)$ unchanged.

Note, however, that our class of dynamics includes examples that are not MRF's (see Example 5 in Section 2.2).

We always assume that the Glauber dynamics is *irreducible*, that is, for any two feasible configurations $\sigma, \tau$, we have $P^t(\sigma, \tau) > 0$ for some finite $t$, and *aperiodic*, that is, for all feasible $\sigma, \tau$, we have $\gcd\{t : P^t(\sigma, \tau) > 0\} = 1$. (Note that aperiodicity is essentially a trivial requirement, and can be enforced by introducing a uniform self-loop probability everywhere.) These two conditions, together with reversibility, ensure that the Glauber dynamics $(X_t)_{t=0}^\infty$ converges to the stationary distribution $\pi$ as $t \to \infty$, for any feasible initial configuration $X_0$. We measure the rate of convergence by the *mixing time*,

$$(1) \qquad \tau_{\mathrm{mix}} := \max_{X_0} \min\left\{t : \|X_t - \pi\| \leq \frac{1}{2e}\right\},$$

where $\|\cdot\|$ denotes total variation distance[5] and (with some abuse of notation) we identify the random variable $X_t$ with its distribution. The constant $\frac{1}{2e}$ is chosen for algebraic convenience only; this choice ensures (see [2]) that $\|X_t - \pi\| \leq \varepsilon$ for all $t \geq \lceil \ln \varepsilon^{-1} \rceil \tau_{\mathrm{mix}}$. In particular, we have

$$(2) \qquad \|X_{k\tau_{\mathrm{mix}}} - \pi\| \leq e^{-k}$$

for all natural numbers $k$.

Finally, we will always assume that the Glauber dynamics satisfies the following *nonredundancy* condition: for every site there are at least two distinct spin values that appear in some feasible configuration. Without such a condition, the system could be artificially padded with an arbitrary number $n'$ of additional "frozen" sites while increasing the mixing time by a factor of only $(n+n')/n$. Equivalently, the mixing time of the continuous-time dynamics (see Section 2.3) would be unchanged.

---

[5]The *total variation distance* between two probability distributions $\mu, \eta$ on a finite set $\Omega$ is defined by $\|\mu - \eta\| := \frac{1}{2} \sum_{z \in \Omega} |\mu(z) - \eta(z)|$.

2.2. *Examples.* We now present a few well-known examples that fit into the above framework. The first four are Markov random fields, so the heat-bath or Metropolis dynamics applies to them (and both are completely specified by the distribution $\pi$). The final example shows that the framework also includes systems that are not MRFs.

1. *The Ising model.* Here $Q = \{\pm 1\}$, and the distribution over configurations is given by

$$\pi(\sigma) \propto \exp\left\{\beta\left(\sum_{\{u,v\} \in E} \sigma(u)\sigma(v) + h \sum_{u \in U} \sigma(u)\right)\right\},$$

where $\beta$ is the inverse temperature and $h$ is the external field. There are no hard constraints. This distribution assigns large weight to configurations in which many neighboring spins are aligned, and many spins have the same sign as $h$. The amount by which such configurations are favored is determined by $\beta$. Plainly this is an MRF as the conditional probability of a configuration on any subset $U$ depends only on the spins on the neighbors of $U$[6]. We note also that this (and the other models we consider) can be extended to the case of arbitrary *boundary conditions*, that is, fixed spins on certain sites of $G$. More precisely, to do this we remove the sites with fixed spins (which are redundant) and incorporate the effects of these spins into the potentials on the neighboring sites.

2. *The hard-core model* (*independent sets*). Here $Q = \{0, 1\}$, and we call sites with spins 0 and 1 "unoccupied" and "occupied" respectively. The feasible configurations $\sigma$ are those in which no two adjacent sites are occupied, that is, they are just the independent sets in $G$; thus there are hard constraints. A feasible configuration with $k$ occupied sites is assigned probability proportional to $\lambda^k$, where the parameter $\lambda > 0$ controls the density of occupation.

3. *Graph colorings* (*the zero-temperature antiferromagnetic Potts model*). Here $Q = \{1, \ldots, q\}$ is thought of as a set of *colors*, and configurations as vertex colorings of $G$. The feasible configurations are proper colorings, that is, those in which no two adjacent vertices receive the same color. The distribution $\pi$ is uniform over proper colorings. Note that the heat-bath (or Metropolis) dynamics may not be irreducible in this case; for instance, there may exist proper colorings that are "frozen," in the sense that no single site can have its color changed (even though other proper colorings do exist). However, it is easy to check that if there are sufficiently many colors, specifically if $q \geq \Delta + 2$ where $\Delta$ is the maximum degree of $G$, then the dynamics is ergodic over all proper colorings.

---

[6]Note that this MRF is completely specified by "potentials" on the sites and edges of $G$, that is, it is a "Gibbs distribution" familiar from statistical physics. [Here the site potentials are $\phi_v(\sigma(v)) = -h\beta\sigma(v)$ for all $v \in V$, and the edge potentials are $\phi_{uv}(\sigma(u), \sigma(v)) = -\beta\sigma(u)\sigma(v)$ for all $\{u, v\} \in E$.] The Hammersley–Clifford theorem (see, e.g., [4]) ensures that such a representation—generalized to include a potential function $\phi_C$ for each *clique C* of $G$—is in fact always possible for an MRF with no hard constraints. If the potentials are allowed to approach infinity, the same holds for MRF's with hard constraints under additional assumptions; see [15].

4. *Constraint satisfaction problems.* Let $\{x_1, \ldots, x_m\}$ be a set of Boolean variables, and let $C_1, \ldots, C_r$ be a set of *constraints*; that is, each $C_i$ specifies, for some subset of the variables, a set of allowed combinations of truth values for those variables. (A canonical example is the Satisfiability problem $k$-SAT, in which each constraint is of the form $z_{i_1} \vee z_{i_2} \vee \cdots \vee z_{i_k}$, where each $z_{i_j}$ denotes either a variable or its negation.) An associated MRF is defined as follows. The graph $G$ has one site for each variable, and an edge between a pair of variables iff they appear together in some constraint. The spin values are $Q = \{\text{True}, \text{False}\}$, and feasible configurations are truth assignments that satisfy all the constraints. The distribution $\pi$ is uniform over feasible configurations. In this example the Glauber dynamics may again fail to be irreducible; this is frequently overcome by considering a "soft" version in which all configurations are feasible but are penalized according to the number of violated constraints they contain.

5. *A non-Markovian example.* It is not hard to construct instances of Glauber dynamics in which the stationary distribution $\pi$ is not an MRF. For a simple example, let $G$ be the cycle on $n$ vertices, and $Q = \{0, 1\}$. Feasible configurations are those in which the set of spin-0 sites form a contiguous block of size $1 \leq \ell \leq n - 1$. The update rule for the spin at any site $v$ is the following: if the spins at the two neighbors of $v$ differ, then flip the spin at $v$, else leave it unchanged. It is easy to check that this dynamics is ergodic and reversible with respect to the uniform distribution on feasible configurations. However, this distribution is clearly not a MRF.

2.3. *Continuous-time dynamics.* It will be convenient in our proofs to work with a continuous-time version of the Glauber dynamics, which we now describe. Given the discrete-time dynamics $(X_t^{\mathcal{D}})_{t=0}^{\infty}$, the continuous-time version $(X_t^{\mathcal{C}})_{t \geq 0}$ is defined as follows:

- there is an independent, rate-1 Poisson clock associated with each site;
- when the clock at $v$ rings, the spin at $v$ is updated as in the discrete-time dynamics.

This process can be viewed equivalently as follows, where $n = |V|$:

- there is a single, rate-$n$ Poisson clock associated with $G$;
- when the clock rings, a site $v$ is chosen u.a.r. and an update is performed at $v$.

It is easy to verify that in distribution $X_t^{\mathcal{C}}$ is equal to $X_{N(nt)}^{\mathcal{D}}$, where $N(nt)$ is the number of rings of a rate-$n$ Poisson clock in time interval $t$. In other words, for every $t \geq 0$,

$$X_t^{\mathcal{C}} \stackrel{\mathrm{d}}{=} \sum_{s=0}^{\infty} e^{-nt} \frac{(nt)^s}{s!} X_s^{\mathcal{D}}. \tag{3}$$

The mixing time for the continuous dynamics is defined by analogy with (1) as

$$\tau_{\mathrm{mix}}^{\mathcal{C}} := \max_{X_0} \inf \left\{ t : \|X_t^{\mathcal{C}} - \pi\| \leq \frac{1}{2e} \right\}. \tag{4}$$

The following observation guarantees that lower bounds on the mixing time translate from the continuous-time to the discrete-time setting. Hence we may work in continuous time without loss of generality.

PROPOSITION 2.1. *For every $t \geq 0$, $\|X_t^{\mathcal{D}} - \pi\| \geq \|X_{2t/n}^{\mathcal{C}} - \pi\| - 2e^{-t}$.*

PROOF. Note first that by (3) and the triangle inequality we may write

$$\|X_{2t/n}^{\mathcal{C}} - \pi\| \leq \sum_{s=0}^{\infty} \mathbf{Pr}(N(2t) = s) \|X_s^{\mathcal{D}} - \pi\|.$$

Splitting the range of summation at $t$, and using the fact that $\|X_s^{\mathcal{D}} - \pi\|$ is a decreasing function of $s$, we find

$$\|X_{2t/n}^{\mathcal{C}} - \pi\| \leq \mathbf{Pr}(N(2t) < t) + \|X_t^{\mathcal{D}} - \pi\|.$$

The desired result follows because, for every $t \geq 0$, $\mathbf{Pr}(N(2t) \leq t) \leq 2e^{-t}$, an elementary property of the Poisson distribution. □

COROLLARY 2.2. *The mixing times for the discrete- and continuous-time dynamics satisfy*

$$\tau_{\mathrm{mix}} \geq \frac{n \tau_{\mathrm{mix}}^{\mathcal{C}}}{6}.$$

PROOF. Set $X_0$ to achieve the maximum in (4). Recall from (2) that

$$\|X_{3\tau_{\mathrm{mix}}}^{\mathcal{D}} - \pi\| \leq \frac{1}{e^3}.$$

Hence, by Proposition 2.1 applied with $t = 3\tau_{\mathrm{mix}}$, and since $\tau_{\mathrm{mix}} \geq 1$,

$$\|X_{6\tau_{\mathrm{mix}}/n}^{\mathcal{C}} - \pi\| \leq \frac{1}{e^3} + 2e^{-3\tau_{\mathrm{mix}}} \leq \frac{3}{e^3} < \frac{1}{2e}$$

which implies $\tau_{\mathrm{mix}}^{\mathcal{C}} \leq 6\tau_{\mathrm{mix}}/n$. □

Corollary 2.2 says that any lower bound on the mixing time of the continuous-time dynamics translates immediately to discrete time, with the loss of only a small constant factor. The translation involves a uniform scaling factor of $n$, arising from the differing clock speeds: while the discrete-time dynamics hits any given site only about once in $n$ steps, the continuous-time version updates the sites at a constant rate.

**3. Basic ingredients.** In this section we introduce two basic tools that will play important roles in our proof. The first bounds the rate at which information can percolate in the Glauber dynamics, and the second expresses a useful monotonicity property of the probability of being in a certain set of configurations. Neither of these is new per se, but our application requires some refinements of them (notably, Lemmas 3.4 and 3.5). We work throughout in continuous time; Section 2.3 shows how to adapt our results to the discrete-time setting.

3.1. *Disagreement percolation.* Suppose $(X_t)$ and $(Y_t)$ are two copies of the Glauber dynamics which agree at time 0 except on some subset $A$ of the sites. Let $A'$ be another subset of sites at distance $d$ from $A$ (so that, in particular, $X_0$ and $Y_0$ agree on $A'$). If $t$ is not too large, then we would expect the distributions of the spin configurations on $A'$ in $X_t$ and $Y_t$ not to differ too much. The extent to which they differ as a function of $t$ is a measure of the rate of information flow in the dynamics, since any difference must be caused by a disagreement percolating from $A$ to $A'$.

We can bound this effect by *coupling* the evolutions of $X_t$ and $Y_t$. Specifically, we will use a *greedy* coupling: we make the two processes use the same Poisson clocks and, whenever the clock at site $v$ rings, we update it to the same value in both processes with the largest possible probability. Thus in particular, if at the time the clock rings $X_t$ and $Y_t$ agree on $v$ and all its neighbors, then this will also hold with probability 1 after the update.

The following lemma bounds the probability that a disagreement percolates from $A$ to $A'$ under this coupling. For each $S \subset V$, recall that $\partial S := \{w \in V \setminus S : \exists u \in S \text{ with } \{u, w\} \in E\}$ denotes the boundary of $S$, and define the "internal boundary" of $S$ by $\delta S := \partial(V \setminus S)$. (Note that $\delta S \subseteq S$.)

LEMMA 3.1. *Let $(X_t)$ and $(Y_t)$ be continuous-time Glauber dynamics on a graph $G = (V, E)$ of maximum degree at most $\Delta$. Suppose $X_0 = Y_0$ at all sites in $V \setminus A$. Let $A' \subset V$ be a set of sites at distance $d > 0$ from $A$. Then the greedy coupling of $(X_t)$ and $(Y_t)$ satisfies*

$$(5) \qquad \mathbf{Pr}(X_t = Y_t \text{ on } A') \geq 1 - \min\{|\delta A|, |\delta A'|\} \left(\frac{et\Delta}{d}\right)^d.$$

*Moreover, the same conclusion holds even if the spin update probabilities of $(X_t)$ and $(Y_t)$ differ at sites in $A$.*

We will often apply Lemma 3.1 in the situation where $A' = B_r(v)$ (the ball of radius $r$ in $G$ centered at $v$), and $A = V \setminus B_{R-1}(v)$ for $R > r$ [so that the evolutions of $(X_t)$ and $(Y_t)$ in the ball $B_{R-1}(v)$ are the same except for possible effects originating from the boundary $S_R(v) := \partial B_{R-1}(v)$]. In this case, the bound in (5) implies

$$\text{(6)} \qquad \mathbf{Pr}(X_t = Y_t \text{ on } A') \geq 1 - \left(\frac{et}{R-r}\right)^{R-r} \Delta^R.$$

The proof of Lemma 3.1 relies on the following simple fact.

OBSERVATION 3.2. *Let $t \geq 0$ and $r \geq 1$. Consider $r$ fully independent Poisson clocks of rate 1. Then the probability, $p$, that there is an increasing sequence of times $0 < t_1 < \cdots < t_r < t$ such that clock $i$ rings at time $t_i$ satisfies*

$$p < \left(\frac{et}{r}\right)^r.$$

PROOF. Since the waiting times between each successive pair of rings of a single clock are independent exponential random variables with mean 1, the event in the lemma has the same probability as the event that a single rate-1 Poisson clock has at least $r$ rings by time $t$. The probability of this event is thus

$$p = e^{-t}\left(\frac{t^r}{r!} + \frac{t^{r+1}}{(r+1)!} + \cdots\right) < \frac{t^r}{r!} < \left(\frac{et}{r}\right)^r.$$

The second inequality here follows because $r! > (\frac{r}{e})^r$ for $r \geq 1$. □

PROOF OF LEMMA 3.1. By the initial conditions and the properties of the greedy coupling, $X_t$ and $Y_t$ can differ at a site in $A'$ only if there exists a "path of disagreement" from $A$ to $A'$. This is a path $v_1, \ldots, v_r$ in $G$ from a site $v_1 \in \delta A$ to a site $v_r \in \delta A'$, together with a sequence of times $t_1 < \cdots < t_r < t$ such that an update is performed at $v_i$ at time $t_i$. Since all such paths have length at least $d$, we may restrict attention to their initial segments of length $d$. For any such initial segment, by Observation 3.2 the probability that the corresponding updates are performed is bounded above by $(\frac{et}{d})^d$. The number of such segments is clearly at most $|\delta A|\Delta^d$. We may also replace $\delta A$ by $\delta A'$ by considering the paths in reverse order. A union bound now completes the proof. Finally, note that the argument did not rely on any properties of the update probabilities for sites in $A$. □

We now derive two further lemmas which use similar ideas but which apply only to dynamics with hard constraints. Call a site *frozen* in a given configuration if its current spin is the only feasible value given the spins on its neighbors. These lemmas bound the rate at which a set of frozen spins can become unfrozen, or vice versa.

LEMMA 3.3. *Let $(X_t)$ be continuous-time Glauber dynamics on a graph $G = (V, E)$ of maximum degree at most $\Delta$. Suppose $X_0$ is frozen except at sites in $A \subset V$. If $A' \subset V$ is a set of sites at distance $d > 0$ from $A$, then*

$$\mathbf{Pr}(X_t = X_0 \text{ on } A') \geq 1 - \min\{|\delta A|, |\delta A'|\}\left(\frac{et\Delta}{d}\right)^d.$$

PROOF. Same as the proof of Lemma 3.1.  □

LEMMA 3.4. *Let $(X_t)$ be continuous-time Glauber dynamics on a graph $G = (V, E)$ of maximum degree at most $\Delta$. Suppose $A, A' \subset V$ are sets of sites at distance $d > 0$ from each other. Then*

$$\mathbf{Pr}((X_t \text{ is frozen on } V \setminus A) \text{ and not } (X_t = X_0 \text{ on } A'))$$
$$\leq \min\{|\delta A|, |\delta A'|\}\left(\frac{et\Delta}{d}\right)^d.$$

PROOF. Label each site with the last time in the interval $[0, t]$ that its spin actually changed. If $X_t$ is frozen on $V \setminus A$, then for any changed site there must be an increasing sequence of such times along a path from that site to $A$. Thus if in addition $X_t$ is not identical to $X_0$ on $A'$, there must be such an increasing sequence of times along a path from $A'$ to $A$. The result now follows in the same way as in the proof of Lemma 3.1.  □

3.2. *Monotonicity properties.* In this section, we will prove the following monotonicity result for Glauber dynamics, which essentially says that, under suitable initial conditions, the probability of a particular site retaining its initial spin after $t$ steps decays slowly as a function of $t$. For this property, we will relax our usual assumptions and require only that the Glauber dynamics is reversible, not necessarily ergodic. This will be the second key ingredient in the proof of Theorem 1.1 in the next section.

LEMMA 3.5. *Let $G = (V, E)$ be a graph, let $v \in V$, and let $Q_v \subseteq Q$. Consider any continuous-time Glauber dynamics on $G$ with spin space $Q$, reversible with respect to a distribution $\pi$ over $\Omega \subseteq Q^V$. Let $\mu$ be the probability, under $\pi$, that the spin at $v$ is in $Q_v$, and suppose $0 < \mu < 1$. Sample the initial configuration $X_0$ according to $\pi$, conditioned on the event that $X_0(v) \in Q_v$. Then, for every $t \geq 0$,*

$$\mathbf{Pr}(X_t(v) \in Q_v) \geq \mu + (1 - \mu)\exp(-t/(1-\mu)).$$

Lemma 3.5 can be derived from a more general monotonicity property of reversible Markov chains. Following Aldous and Fill ([3], Chapter 3, Section 4), call a function $f$ *completely monotone decreasing* (*CMD*) if it can be written in the form

$$f(t) = \sum_{i=1}^{m} \alpha_i \exp(-\lambda_i t), \tag{7}$$

where all the coefficients $\alpha_i, \lambda_i$ are nonnegative.

LEMMA 3.6. *Let $(X_t)$ be a continuous-time Markov chain on finite state space $\Omega$, and let $(X_t)$ be reversible with respect to a distribution $\pi$. Let $\Psi \subset \Omega$. If the distribution of $X_0$ is supported on $\Psi$, and proportional to $\pi$ on $\Psi$, then $\mathbf{Pr}(X_t \in \Psi)$ is CMD. Indeed,*

$$\mathbf{Pr}(X_t \in \Psi) = \sum_{i=1}^{|\Omega|} \alpha_i \exp(-\lambda_i t),$$

*where $\sum_i \alpha_i = 1$, $\alpha_1 = \pi(\Psi)$, $\lambda_1 = 0$ and all $\alpha_i, \lambda_i \geq 0$.*

As mentioned earlier, this (and indeed more general) monotonicity properties for both continuous- and discrete-time Markov chains are well known (see, e.g., [3], Chapter 3, Section 4). We include a proof for completeness.

PROOF OF LEMMA 3.6. As stipulated in the lemma, we let $(X_t)$ have the initial distribution $\mu_0(i) = \pi(i)/\pi(\Psi)$ for $i \in \Psi$ (and 0 elsewhere). Writing $P_t(i,j)$ for the time-$t$ transition probability from $i$ to $j$, our main goal is to show that the quantity

$$\mathbf{Pr}(X_t \in \Psi) = \sum_{i \in \Omega} \sum_{j \in \Psi} \mu_0(i) P_t(i,j) = \pi(\Psi)^{-1} \sum_{i \in \Psi} \sum_{j \in \Psi} \pi(i) P_t(i,j) \tag{8}$$

is CMD.

By the standard spectral representation for reversible Markov chains, we may write

$$P_t(i,j) = \sqrt{\frac{\pi(j)}{\pi(i)}} \sum_{k=1}^{|\Omega|} \exp(-\lambda_k t) u_{ik} u_{jk},$$

where $0 = \lambda_1 \leq \lambda_2 \leq \cdots \leq \lambda_{|\Omega|}$ are the eigenvalues of the transition kernel, and $(u_{ij})_{i,j=1}^{|\Omega|}$ is an orthonormal matrix. Moreover, we may take $u_{i1} = \sqrt{\pi(i)}$, since $\pi$ is a stationary vector for the Markov chain. Substitution into (8) gives

$$\mathbf{Pr}(X_t \in \Psi) = \pi(\Psi)^{-1} \sum_k \exp(-\lambda_k t) \sum_{i \in \Psi} \sum_{j \in \Psi} \sqrt{\pi(i)\pi(j)} u_{ik} u_{jk}$$

$$(9) \qquad = \pi(\Psi)^{-1} \sum_k \exp(-\lambda_k t) \left( \sum_{i \in \Psi} \sqrt{\pi(i)} \, u_{ik} \right)^2.$$

This is clearly of the form (7), so the function is CMD. Moreover, $\alpha_1$ satisfies

$$\alpha_1 = \pi(\Psi)^{-1} \left( \sum_{i \in \Psi} \pi(i) \right)^2 = \pi(\Psi),$$

and by (9) the sum of all coefficients is

$$\sum_{k=1}^{|\Omega|} \alpha_k = \pi(\Psi)^{-1} \sum_k \left( \sum_{i \in \Psi} \sqrt{\pi(i)} u_{ik} \right)^2 = \mathbf{Pr}(X_0 \in \Psi) = 1. \qquad \square$$

We now prove the main result of this section.

PROOF OF LEMMA 3.5. Let $(X_t)$ be as described in the lemma, and define $f(t) := \mathbf{Pr}(X_t(v) \in Q_v)$. By Lemma 3.6, we may write

$$f(t) = \sum_{i=1}^{|\Omega|} \alpha_i \exp(-\lambda_i t),$$

where $\sum_{i=1}^{|\Omega|} \alpha_i = 1$, $\alpha_1 = \mu$, $\lambda_1 = 0$ and all $\alpha_i, \lambda_i \geq 0$.

Now, let us sample $J$ from $1, \ldots, |\Omega|$ according to the distribution $(\alpha_j)$. Note that, in this case,

$$f(t) = \mathbf{E}(\exp(-\lambda_J t)) = \mu + (1-\mu)\mathbf{E}(\exp(-\lambda_J t) \mid J \geq 2).$$

Since the function $\exp(-x)$ is convex, we may apply Jensen's inequality to the conditional probability on the right-hand side to obtain

$$f(t) \geq \mu + (1-\mu)\exp(-t\mathbf{E}(\lambda_J \mid J \geq 2)).$$

Observe that

$$\mathbf{E}(\lambda_J \mid J \geq 2) = \frac{1}{1-\mu} \sum_{j=2}^{|\Omega|} \alpha_j \lambda_j = \frac{-f'(0)}{1-\mu}.$$

Moreover, for every $s \geq 0$,

$$f(s) \geq \mathbf{Pr}(\text{clock at } v \text{ does not ring by time } s) = e^{-s},$$

and hence $f'(0) \geq -1$, which completes the proof. $\square$

**4. Main result.** We are now ready to state and prove the continuous-time version of our main result.

THEOREM 4.1. *Let $\Delta \geq 2$ be a fixed positive integer, and let $G$ be any graph on $n$ vertices with maximum degree at most $\Delta$. Any continuous-time Glauber dynamics on $G$ has mixing time $\Omega(\log n)$, where the constant in the $\Omega(\cdot)$ depends only on $\Delta$.*

By Corollary 2.2, Theorem 4.1 immediately implies that the discrete-time Glauber dynamics has mixing time at least $\Omega(n \log n)$, thus proving our main result (Theorem 1.1) claimed in the Introduction.

We prove Theorem 4.1 in two steps. First, in Section 4.1, we consider the conceptually simpler case in which there are no hard constraints. Then, in Section 4.2, we show how to extend the proof to handle hard constraints.

4.1. *Proof for the case of no hard constraints.* Let $T = \frac{\ln n}{8e\Delta \ln \Delta}$. We will show that, for a suitable initial distribution $X_0$, the variation distance between $X_T$ and the stationary distribution is large and hence $\tau_{\text{mix}}^{\mathcal{C}} \geq T$. In the remainder of the proof, we will assume that $\Delta \geq 2$ is fixed and $n$ is sufficiently large.

Let $R = \lceil \frac{\ln n}{4 \ln \Delta} \rceil$. Choose a set of $\lceil \frac{n}{\Delta^{2R}} \rceil$ sites $C \subset V$ whose pairwise distances are all at least $2R$. (Such a set $C$ can always be found by repeatedly picking a site, adding it to $C$, and removing it from $V$, together with all sites at distance $\leq 2R - 1$.) We call the elements of $C$ *centers*, since we will shortly be considering a restriction of the dynamics to balls of radius $R-1$ centered at elements of $C$.

For each site $v \in C$, let $Q_v$ be a nonempty proper subset of the spins. (For instance, they could all be the singleton $\{1\}$.) For every configuration $X$, let $f(X)$ denote the fraction of centers $v$ such that $X(v) \in Q_v$. We will specify a distribution for $X_0$ and a threshold $\widehat{\mu}$, such that

$$\left| \mathbf{Pr}(f(X_T) \geq \widehat{\mu}) - \Pr_{X \sim \pi}(f(X) \geq \widehat{\mu}) \right| > \frac{1}{2e},$$

where $\pi$ is the stationary distribution. This implies that the mixing time is greater than $T$.

Let $U$ denote the set of sites at distance $\geq R$ from $C$, and let $\sigma_U$ be an arbitrary assignment of spins to all sites in $U$. The distribution of $X_0$ will be concentrated on configurations which agree with $\sigma_U$, and satisfy either $(\forall v \in C)(X_0(v) \in Q_v)$, or $(\forall v \in C)(X_0(v) \notin Q_v)$ (we will specify which in the next paragraph). Among configurations satisfying these constraints, the distribution of $X_0$ will be proportional to $\pi$.

Let $\mu$ denote the conditional expectation of $f(Y)$ when $Y$ is drawn from the stationary distribution $\pi$, *conditioned on agreeing with $\sigma_U$*. Let $\varepsilon = 1/(4\exp(2T))$, and define $\widehat{\mu} \in [0,1]$ by

(10) $$\widehat{\mu} = \begin{cases} \mu - \varepsilon, & \text{if } \mu > 1/2, \\ 1/2, & \text{if } \mu = 1/2, \\ \mu + \varepsilon, & \text{if } \mu < 1/2. \end{cases}$$

Let $X \sim \pi$. If $\mathbf{Pr}(f(X) \geq \widehat{\mu}) \leq 1/2$, then we require that $(\forall v \in C)(X_0(v) \in Q_v)$. Otherwise, we require that $(\forall v \in C)(X_0(v) \notin Q_v)$. Replacing $Q_v$ by its complement if necessary, we may assume without loss of generality that $\mathbf{Pr}(f(X) \geq \widehat{\mu}) \leq 1/2$ and $(\forall v \in C)(X_0(v) \in Q_v)$.

With this definition of $X_0$, it suffices to show that

(11) $$\mathbf{Pr}(f(X_T) < \widehat{\mu}) < \frac{1}{2} - \frac{1}{2e}.$$

To this end, we introduce a second chain, $Y_0, \ldots, Y_T$, which evolves according to the Glauber dynamics except that only the spins on the balls $\bigcup_{v \in C} B_{R-1}(v)$ are updated. (The sites outside these balls retain their initial configurations and ignore their clocks.) Note that this chain decomposes into fully independent processes on each of the balls $B_{R-1}(v)$. For the initial distribution, we take $Y_0 = X_0$.

We will first show that, under the greedy coupling of $(X_t)$ with $(Y_t)$, the expected number of centers $v \in C$ at which $X_T$ and $Y_T$ differ is small, and hence that

(12) $$\mathbf{Pr}(f(X_T) \leq f(Y_T) - \varepsilon/2) \leq 1/4.$$

Then we will take advantage of the independence of the spins $Y_T(v), v \in C$, to show that $f(Y_T)$ is concentrated around its expectation and hence also that

(13) $$\mathbf{Pr}(f(Y_T) < \widehat{\mu} + \varepsilon/2) < 1/4 - 1/(2e).$$

Combining inequalities (12) and (13) yields (11) and thus will complete the proof.

Fix $v \in C$. Since $X_0 = Y_0$, and the update rules for the two chains are the same on $B_{R-1}(v)$, it follows from Lemma 3.1 and equation (6) (with $r = 0$) that, under the greedy coupling,

$$\mathbf{Pr}(X_T(v) \neq Y_T(v)) \leq \left(\frac{eT\Delta}{R}\right)^R \leq \varepsilon/8.$$

The second inequality here follows by plugging in the values for $T$, $R$ and $\varepsilon$. By linearity of expectation, it follows that

$$\mathbf{E}(|\{v \in C : X_T(v) \neq Y_T(v)\}|) \leq |C|\varepsilon/8,$$

for all sufficiently large $n$. Applying Markov's inequality to the random variable $|\{v \in C : X_T(v) \neq Y_T(v)\}|$ now yields our first desired inequality (12).

Again fix $v \in C$, and let us restrict our attention to the ball $B_{R-1}(v)$. Let $\pi_{\sigma_U}$ denote the stationary distribution $\pi$ on $B_{R-1}(v)$, conditioned on agreeing with $\sigma_U$. Let $\mu_v = \mathbf{Pr}(Y(v) \in Q_v)$, where $Y$ is sampled from $\pi_{\sigma_U}$. Note that the dynamics $(Y_t)$ is reversible with respect to this distribution (although it need not be ergodic, because of the fixed boundary). Note also that the initial distribution $Y_0$ is $\pi_{\sigma_U}$ conditioned on $Y_0(v) \in Q_v$. Hence by Lemma 3.5 applied to $(Y_t)$ on the ball $B_{R-1}(v)$, we have

$$\mathbf{Pr}(Y_T(v) \in Q_v) \geq \mu_v + (1 - \mu_v) \exp(-T/(1 - \mu_v)).$$

There is a subtle technical point in the preceding argument: since we are working in continuous time, and because the boundary of $B_{R-1}(v)$ is fixed, the projection of $(Y_t)$ onto $B_{R-1}(v)$ is indeed a Glauber dynamics and so Lemma 3.5 applies.

By linearity of expectation, $\mu = \frac{1}{|C|} \sum_{v \in C} \mu_v$. Also by linearity of expectation, we have

$$\mathbf{E}(f(Y_T)) = \frac{1}{|C|} \sum_{v \in C} \mathbf{Pr}(Y_T(v) \in Q_v)$$

$$\geq \mu + \frac{1}{|C|} \sum_{v \in C} (1 - \mu_v) \exp(-T/(1 - \mu_v)).$$

Since the function $x \exp(-t/x)$ is convex on $(0, 1)$ as a function of $x$, for fixed $t > 0$, it follows by Jensen's inequality that

(14) $$\mathbf{E}(f(Y_T)) \geq \mu + (1 - \mu) \exp(-T/(1 - \mu)).$$

We now claim that the right-hand side of (14) is at least $\widehat{\mu} + \varepsilon$. To see this, note by the definition of $\widehat{\mu}$ that, when $\mu > 1/2$, we have $\mu = \widehat{\mu} + \varepsilon$. On the other hand, when $\mu \leq 1/2$, we have $(1 - \mu) \exp(-T/(1 - \mu)) \geq 1/(2 \exp(2T)) = 2\varepsilon$ and $\mu \geq \widehat{\mu} - \varepsilon$. We conclude from (14) that $\mathbf{E}(f(Y_T)) \geq \widehat{\mu} + \varepsilon$.

Since $f(Y_T)$ is the average of $|C|$ independent random variables taking values in $\{0, 1\}$, Chernoff's bound yields

$$\mathbf{Pr}(f(Y_T) < \widehat{\mu} + \varepsilon/2) \leq \mathbf{Pr}(f(Y_T) < \mathbf{E}(f(Y_T)) - \varepsilon/2)$$

$$\leq \exp(-|C|\varepsilon^2/8).$$

Writing $|C|$ and $\varepsilon$ explicitly in terms of $n$ shows that the right-hand side is asymptotically 0 as $n \to \infty$, establishing (13).

4.2. *Proof for the general case.* The proof of the general version of Theorem 4.1 will rely on the following fact about spin systems with hard constraints.

LEMMA 4.2. *Let $v \in V$ and let $R \geq 4$ be an integer. Let $(X_t)$ and $(Y_t)$ be continuous-time Glauber dynamics such that the initial configuration $X_0 \in \Omega$ is frozen at all sites in $B_R(v)$, and the initial configuration $Y_0$ differs from $X_0$ at $v$. If $T \leq R/(5e^2 \Delta \ln \Delta)$, then the distributions of $X_T$ and $Y_T$ are at total variation distance at least*

$$1 - 2\exp(-R/(3\ln \Delta)).$$

PROOF. Let $0 < r < R$ be another positive integer. By Lemma 3.3 applied to $(X_t)$ [setting $A = V \setminus B_R(v)$ and $A' = B_r(v)$], the probability that $X_T$ agrees with $X_0$ on $B_r(v)$ is at least $1 - |\delta A'|(eT\Delta/(R-r+1))^{R-r+1} \geq 1 - \Delta(\Delta-1)^{r-1}(eT\Delta/(R-r+1))^{R-r+1}$. On the other hand, by applying Lemma 3.4 to $(Y_t)$ [setting $A = V \setminus B_r(v)$ and $A' = \{v\}$], the probability that $Y_T$ agrees with $X_0$ on $B_r(v)$ is at most $(eT\Delta/(r+1))^{r+1}$. By the triangle inequality, the distributions of $X_T$ and $Y_T$ must have total variation distance at least

$$1 - \Delta^r(eT\Delta/(R-r+1))^{R-r+1} - (eT\Delta/(r+1))^{r+1}.$$

Letting $r = \lfloor R/3\ln(\Delta) \rfloor < R/2$, the result follows by basic algebra. □

We are now ready to present the proof of Theorem 4.1 in the presence of hard constraints. The broad outline of the proof follows that of the previous one where there were no hard constraints; we will focus on the points where differences occur. In particular, the selection of the set of centers $C$, the spin sets $Q_v$, and the initial assignment $\sigma_U$ to the sites at distance $> R$ from $C$, become slightly trickier.

Let $R = \lceil \frac{\ln n}{4 \ln \Delta} \rceil$ as in the previous proof.

*Case* 1. Suppose there exists a set $C \subset V$ of size at least $\lceil n/\Delta^{3R} \rceil = \Theta(n^{1/4})$ at pairwise distance at least $2R$, together with an initial assignment $\sigma_U$ to all the sites $U$ at distance $\geq R$ from $C$, such that, conditioned on the configuration agreeing with $\sigma_U$, every $v \in C$ has at least two feasible spins available.

In this case, for each $v \in C$, define $Q_v$ to be a proper nonempty subset of the feasible spins for $v$, conditioned on agreeing with $\sigma_U$. Then the rest of the proof goes through as in the previous subsection, but with slightly worse constants.

*Case* 2. No sufficiently large set $C$ as above exists.

In this case, let $C$ and $\sigma_U$ be as above, with $|C| < n/\Delta^{3R}$ and $C$ maximal. The upper bound on $|C|$ implies that there must exist at least one site $v$ at distance $\geq 3R$ from $C$. So, in particular, every site in $B_R(v)$ is at distance $\geq 2R$ from $C$. It follows that all of $B_R(v)$ must be frozen under $\sigma_U$; otherwise, we could add any one of its nonfrozen sites to $C$, contradicting maximality.

To complete the proof, we argue differently from the previous case: we exhibit two initial configurations $X_0$ and $Y_0$ such that when both evolve according to the Glauber dynamics for time $T \leq \gamma \log n$ (for sufficiently small $\gamma > 0$), $X_T$ and $Y_T$ are still at large variation distance. Hence the mixing time must be at least $T$. Let $X_0$ be any configuration extending $\sigma_U$, and let $Y_0$ be any configuration which disagrees with $X_0$ at $v$. (Such a $Y_0$ must exist by our nonredundancy condition.) By Lemma 4.2, the mixing time is $\Omega(R/(\Delta \log \Delta)) = \Omega(\log n/(\Delta \log^2 \Delta))$.

This completes the proof of Theorem 4.1 in the general case.

REMARK 4.3. In both cases in the above proof, the dependence of our lower bound on $\Delta$ is $\Omega(1/\Delta \operatorname{polylog} \Delta)$. In the next section, we will prove an upper bound whose dependence on $\Delta$ is $O(1/\log \Delta)$. It would be interesting to close this gap and establish a tight dependence on $\Delta$.

**5. Unbounded degrees and other variations.** In this final section we discuss various extensions of our main theorem, and also show that certain other extensions (notably, removing the assumption of bounded degree) are not possible.

5.1. *Graphs with unbounded degree.* We have shown a lower bound of $\Omega(n \log n)$ on the mixing time of Glauber dynamics on any family of graphs of bounded degree. More precisely, for any fixed $\Delta$, we have shown that the Glauber dynamics on graphs of maximum degree at most $\Delta$ has mixing time at least $C_\Delta n \ln n$ for some positive constant $C_\Delta$. Is the restriction to bounded-degree graphs necessary?

We first give a simple example to show that it is. Let $G = K_n$, and consider the hard-core model on $G$, that is, the spin space is $Q = \{0, 1\}$, and the feasible configurations are those in which at most one vertex has spin 1 (corresponding to the independent sets of $G$). The distribution $\pi$ assigns probability $\frac{1}{2}$ to the all-0 configuration, and $\frac{1}{2n}$ to each other configuration. (Thus the activity parameter in the hard-core model is $\lambda = \frac{1}{n}$.) Consider the following Metropolis Glauber dynamics, which at each step picks a vertex $v$ u.a.r. If $v$ has spin 1, then its spin is flipped to 0 with probability $\frac{1}{2}$, while if $v$ and all other vertices have spin 0, then the spin of $v$ is flipped to 1 with probability $\frac{1}{2n}$. Plainly the mixing time of this Glauber dynamics (in discrete time) is $O(n)$.

The following more interesting example demonstrates a trade-off between maximum degree and mixing time. Let $G_0 = (V, E_0)$ be any graph of maximum degree $d \leq \sqrt{\Delta}$. Let $G = (V, E)$, where $E$ consists of all pairs of vertices whose distance in $G_0$ is 1 or 2 (note that the maximum degree of $G$ is at most $d^2 \leq \Delta$). We will construct a Glauber dynamics on $G$ which mimics any reversible random walk on $G_0$. To this end, for every edge $\{u,v\} \in E_0$ we let $A(u,v) > 0$ be the transition probability from $u$ to $v$ in such a random walk, $A(v,v) \geq 0$ be the probability of a self-loop at $v$, and $\pi_0$ be the stationary distribution. For convenience, in what follows we augment $E_0$ to include the self-loop $\{v,v\}$ whenever $A(v,v) > 0$.

Our Glauber dynamics will have spins $Q = \{0,1,2\}$. The feasible configurations $\sigma: V \to Q$ will all satisfy $\sum_{v \in V} \sigma(v) \in \{1,2\}$. Moreover, when this sum is 2 the vertices with nonzero spin will be the endpoints of an edge $\{u,v\}$ (possibly a self-loop) in $E_0$. Thus we may identify feasible configurations with the set $V \cup E_0$.

We now describe the update rule for the Glauber dynamics. Let $\sigma$ be the current configuration, and let $v$ be the site selected for updating. The new configuration $\sigma'$ is determined as follows (where $\sim$ denotes adjacency in $G_0$):

1. Case 1: $\sigma(v) = 0$.

    (i) If there exists $w \sim v$, $\sigma(w) = 1$, and for all $z \sim w$, $\sigma(z) = 0$, then: with probability $A(w,v)$, set $\sigma'(v) = 1$.
    (ii) Otherwise, set $\sigma'(v) = 0$.

2. Case 2: $\sigma(v) = 1$.

    (i) If there exists $w \sim v$, $\sigma(w) = 1$, then set $\sigma'(v) \in \{0,1\}$ according to the result of a fair coin toss.
    (ii) Otherwise, with probability $A(v,v)$, set $\sigma'(v) = 2$.

3. Case 3: $\sigma(v) = 2$. Set $\sigma'(v) = 1$.

Since the update rule only examines vertices at distance $\leq 2$ from $v$ in $G_0$, this is a local dynamics on the graph $G$.

Informally, the Glauber dynamics may be described as follows. As indicated above, we may think of the configurations of the dynamics as either *vertex states*, $\sigma_v$ for some $v \in V$, or *edge states*, $\sigma_{vw}$ for some $\{v,w\} \in E_0$. Call a transition of the Glauber dynamics *successful* if it causes a change in the configuration. Each successful transition causes a change from a vertex state to an edge state, or vice versa. Starting from a vertex state $\sigma_v$, the first successful transition moves to an edge $\sigma_{vw}$ (where $w = v$ is possible); the next successful transition undoes this move with probability $\frac{1}{2}$, and otherwise moves to the vertex state $\sigma_w$, thus completing the simulation of a single move of the original random walk.

It is easy to check that the stationary distribution of the dynamics is $\pi = \frac{1}{2}(\pi^V + \pi^E)$, where $\pi^V$ is concentrated on vertex states and satisfies $\pi^V(\sigma_v) = \pi_0(v)$ for all $v \in V$, and $\pi^E$ is concentrated on edge states and satisfies $\pi^E(\sigma_{vw}) = \pi_0(v)A(v,w) + \pi_0(w)A(w,v) = 2\pi_0(v)A(v,w)$ for $\{v,w\} \in E_0$ with $v \neq w$, $\pi^E(\sigma_{vv}) = \pi_0(v)A(v,v)$ for self-loops. Since $A$ is reversible this dynamics is also reversible. Thus we have a valid Glauber dynamics.

LEMMA 5.1. *The above Glauber dynamics has mixing time $O(n\tau_A)$, where $\tau_A$ is the mixing time of the underlying random walk defined by $A$ and $n = |V|$.*

PROOF (Sketch). Let $(X_{t_0}, X_{t_1}, \ldots, X_{t_k}, \ldots)$ denote the Glauber dynamics observed after each successful transition, with $t_0 = 0$ and $X_{t_0}$ assumed to be a vertex state. Then, from the informal description given earlier, $(X_{t_{2k}})$ is a (reversible) Markov chain on the vertex states, with transition matrix $\frac{1}{2}(I + A)$, where $I$ is the $n \times n$ identity matrix. Thus, denoting the distribution of $X_{t_k}$ by $\nu_k$, if we take $k \geq C\tau_A$ for a sufficiently large constant $C$, then if $k$ is even we have $\nu_k = \pi^V + \varepsilon_k$ where $\|\varepsilon_k\|_1 \leq \varepsilon$ for any desired positive constant $\varepsilon$. From the definition of the dynamics, this then easily implies $\nu_k = \pi^E + \varepsilon_k$ for $k \geq C\tau_A$, $k$ odd, as well, with $\|\varepsilon_k\|_1$ bounded similarly.

Now let $(X_0, X_1, \ldots, X_t, \ldots)$ denote the full Glauber dynamics, and let $\mu_t$ denote the distribution of $X_t$. Then $\mu_t$ satisfies $\mu_t = \sum_k p_k^{(t)} \nu_k$, where $p_k^{(t)}$ is the probability that exactly $k$ successful transitions occur within the first $t$ steps of the dynamics. Note that the probability of any given transition being successful is exactly $\frac{1}{n}$. Thus if we take $t = C'n\tau_A$ for a suitable constant $C'$, we may ensure that $\sum_{k < C\tau_A} p_k^{(t)} \leq \varepsilon$ and $|\sum_{k \geq C\tau_A,\ k \text{ even}} p_k^{(t)} - \frac{1}{2}| \leq \varepsilon$. This in turn implies

$$\mu_t = \sum_{k < C\tau_A} p_k^{(t)} \nu_k + \sum_{k \geq C\tau_A} p_k^{(t)} \nu_k = \tfrac{1}{2}(\pi^V + \pi^E) + \gamma_t,$$

where $\|\gamma_t\|_1 \leq 3\varepsilon$. □

Our motivation for the general construction above is the following special case, which illustrates a trade-off between maximum degree and mixing time. This is the result we claimed in Theorem 1.2, which we restate here for convenience.

THEOREM 5.2. *For each $n$, let $\Delta(n)$ be any natural number satisfying $2 \leq \Delta(n) < n$. Then there exists a family of graphs $G = G(n)$, where $G(n)$ has $n$ vertices and maximum degree $\Delta(n)$, and an associated Glauber dynamics on $G(n)$ with mixing time $O(n \log n / \log \Delta(n))$.*

PROOF. Let $\Delta = \Delta(n)$ and set $d = \lfloor \sqrt{\Delta} \rfloor$. Without loss of generality we may assume $d \geq 3$, as there are many examples of Glauber dynamics with mixing time $O(n \log n)$ on arbitrary graphs. Let $G_0$ be a complete $(d-1)$-ary tree with $n$ vertices; note that the height of $G_0$ is $h = \Theta(\log n / \log d)$. Let $A$ be the biased random walk on $G_0$ which has, at each step, probability $2/3$ of moving to a parent (except at the root), and $1/3$ of moving to a random child (except at a leaf). It is easily seen (by projecting onto a one-dimensional process that walks between the levels of $G_0$) that the mixing time of this random walk is $\Theta(h) = \Theta(\log n / \log d)$. Therefore, by Lemma 5.1, the Glauber dynamics constructed as above from this random walk has mixing time $O(n \log n / \log d) = O(n \log n / \log \Delta)$ and degree $d^2 \leq \Delta$, as required. (The maximum degree can always be increased to exactly $\Delta$ by adding edges as necessary.) □

REMARK. It is in fact possible, with a little more work, to massage the above tree example into one in which the Glauber dynamics is based on a Markov random field.

5.2. *Nonuniform update probabilities.* In our main theorem, we assumed that the site to be updated was chosen uniformly at random (or equivalently, in continuous time, all sites are updated at rate 1). How essential is this assumption?

First, it is not hard to check that our proof still goes through more or less unchanged under the weaker assumption that no site has a probability larger than $O(1/n)$ of being chosen. On the other hand, if we allow arbitrary site selection probabilities, the theorem no longer holds. To see this, consider again the tree example from Section 5.1, with, for example, $d = 3$, but now suppose that vertex $v$ is chosen with probability proportional to $(d-1)^{k/2}$, where $k$ is the height of $v$ in $G_0$. Thus the probability of choosing $v$ is $c(d-1)^{k/2-h}$, where $c = (1-(d-1)^{-1/2})/(1-(d-1)^{-(h+1)/2}) > 1/4$. The updating rule for $v$ is as before. Now by an analysis similar to that of Section 5.1, the expected time to reach the root from a leaf is approximately proportional to the sum of the waiting times to hit vertices along a leaf-root path, which is at most $\sum_{k=0}^{h} c^{-1}(d-1)^{h-k/2} = O(n)$. This linear mixing time is significantly less than our $\Omega(n \log n)$ lower bound for (near-)uniform site selection probabilities.

5.3. *Block dynamics.* An alternative to Glauber dynamics frequently encountered in the literature is so-called "block dynamics," which, at each step, updates all the spins in a randomly chosen "block" of sites. (For instance, a block might be the ball of radius 1 around a randomly chosen site.) In many cases it is easier to establish rapid mixing for block dynamics (even with quite small blocks) than for single-site dynamics (see, e.g., [1, 11, 16]). As with the Glauber dynamics, we require that the block dynamics be irreducible and aperiodic, reversible with respect to $\pi$, and local in the sense that the updates inside the block depend only on the current spins on the sites in the block and their neighbors.

Provided again that $G$ has bounded degree, the blocks are connected and of bounded size, and no site is updated at a rate greater than $O(1)$ (in continuous time), our proof of Theorem 4.1 can easily be seen to apply in this context as well, giving a lower bound of $\Omega(n \log n)$ on the mixing time of block dynamics.

**Acknowledgments.** We thank Laci Babai, Mark Jerrum and Yuval Peres for enlightening discussions, and Mark Jerrum for sharing a preliminary version of [9]. We would also like to thank Jaikumar Radhakrishnan, as well as the anonymous referees, for helpful suggestions.

An extended abstract of this paper appeared in *Proceedings of the 46th Annual IEEE Symposium on Foundations of Computer Science* (FOCS 2005), pages 511–520.

TOYOTA TECHNOLOGICAL INSTITUTE
1427 E. 60TH ST.
CHICAGO, ILLINOIS 60637
USA
E-MAIL: hayest@tti-c.org

COMPUTER SCIENCE DIVISION
UNIVERSITY OF CALIFORNIA BERKELEY
BERKELEY, CALIFORNIA 94720
USA
E-MAIL: sinclair@cs.berkeley.edu